\documentclass[twoside, 11pt, leqno]{amsart}
 \usepackage[latin1]{inputenc}
 \usepackage[T1]{fontenc}
 \usepackage{amsmath, amssymb, amsfonts, amsthm, amscd}
 \usepackage[left=110pt, right=110pt, bottom=110pt, top=100pt, footskip=22pt]{geometry}
 \usepackage{mathrsfs}
 \usepackage[all]{xy}
 \usepackage{array, booktabs, ragged2e}
 \usepackage{slashed}

\theoremstyle{plain}
\newtheorem{theorem}{Theorem}
\newtheorem{corollary}[equation]{Corollary}
\newtheorem{lemma}[equation]{Lemma}
\newtheorem{proposition}[equation]{Proposition}
\newtheorem{conjecture}[equation]{Conjecture}
\newtheorem{definition}[equation]{Definition}

\newtheorem{remark}[equation]{Remark}

\newtheorem{geoassumptions}[equation]{Geometric assumptions}
\newtheorem{analassump}[equation]{Analytic assumptions}

\numberwithin{equation}{section}
\numberwithin{theorem}{section}

\DeclareMathOperator{\cvf}{cvf}
\DeclareMathOperator{\ran}{ran}
\DeclareMathOperator{\tr}{tr}
\DeclareMathOperator{\TR}{TR}
\DeclareMathOperator{\Vol}{Vol}
\DeclareMathOperator{\Metr}{Metr}

\DeclareMathOperator{\ord}{ord}
\DeclareMathOperator{\Ind}{Ind}
\DeclareMathOperator{\Homo}{Hom}
\DeclareMathOperator{\Irr}{Irr}
\DeclareMathOperator{\Proj}{Proj}

\DeclareMathOperator{\Scal}{Scal}
\DeclareMathOperator{\Ric}{Ric}
\DeclareMathOperator{\sign}{sign}
\DeclareMathOperator{\Hess}{Hess}

\DeclareMathOperator{\Id}{Id}
\DeclareMathOperator{\diver}{div}
\DeclareMathOperator{\conf}{conf}
\DeclareMathOperator{\diff}{diff}
\DeclareMathOperator{\ctran}{Conf}
\DeclareMathOperator{\Diffeo}{Diff}
\DeclareMathOperator{\confdiff}{(conf+diff)_{g_0}}
\DeclareMathOperator{\confdiffperp}{(conf+diff)_{g_0}^\perp}

\newcommand{\Proof}{\begin{proof}}
\newcommand{\Endproof}{\end{proof}}
\newcommand{\fracsm}[2]{\begin{matrix}\frac{#1}{#2}\end{matrix}}
\newcommand{\prefix}[3]{\vphantom{#3}#1#2#3}
\newcommand{\beq}{\begin{equation}}
\newcommand{\eeq}{\end{equation}}

\newcommand{\Naturals}{\mathbb{N}}
\newcommand{\Reals}{\mathbb{R}}
\newcommand{\Complex}{\mathbb{C}}

\newcommand{\Real}{\mathrm{Re}\,}

\newcommand{\SO}[1]{\mathrm{SO}(#1)}
\newcommand{\SU}[1]{\mathrm{SU}(#1)}
\newcommand{\SOzero}[1]{\mathrm{SO}_0(#1)}
\newcommand{\Dirac}{D}

\newcommand{\Diff}[1]{\mathrm{Diff}(#1)}
\newcommand{\DiffO}[1]{\mathrm{Diff}_0(#1)}
\newcommand{\Cinf}[1]{C^\infty(#1)}          
\newcommand{\CinfStwoM}{C^\infty(S^2TM)}

\newcommand{\Endo}[1]{\mathrm{End}(#1)} 
\newcommand{\DIP}[2]{\langle\!\langle #1,#2\rangle\!\rangle}

\newcommand{\Lie}{L}
\newcommand{\TFSM}{C^\infty(S_0^2TM)}
\newcommand{\TFSSn}{C^\infty(S_0^2TS^n)}

\makeindex

\title[Rigidity of conformal functionals]{Rigidity of conformal functionals on spheres}
\author{Niels Martin M\o{}ller}
\author{Bent \O{}rsted}
\address{Niels Martin M\o{}ller, Department of Mathematics, Massachusetts Institute of Technology, Cambridge, MA 02139, USA}
\email{moller@math.mit.edu} 
\address{Bent \O{}rsted, Department of Mathematical Sciences, University of Aarhus, DK-8000 Aarhus C, Denmark}
\email{orsted@imf.au.dk} 
\thanks{The authors would like to thank Alice Chang, Andreas Juhl, Robin Graham, Kate Okikiolu, Peter Sarnak, Paul Yang and others, for many helpful discussions concerning the results in this paper. The first author was supported by Department of Mathematical Sciences, University of Aarhus, Denmark, and partly by the Elite Research Travel Grant 2007, from the Danish Ministry of Science, Technology and Innovation.}

\date{February 2009}

\begin{document}

\begin{abstract}
In this paper we investigate the nature of stationary points of
functionals on the space of Riemannian metrics on a smooth compact
manifold. Special cases are spectral invariants associated with
Laplace or Dirac operators such as functional determinants, and the total $Q$-curvature. When
the functional is invariant under conformal changes of the metric,
and the manifold is the standard $n$-sphere, we apply methods from
representation theory to give a universal form of the Hessian of
the functional at a stationary point. This reveals a very strong rigidity in the local structure of any such functional. As a corollary this gives a new proof of the results of K. Okikiolu (Ann. Math., 2001) on local maxima and minima for the determinant of the conformal Laplacian, and we obtain results of the same type in general examples.
\end{abstract}

\maketitle
\section{Introduction}\label{sec:introduction}
In recent years there has been much progress in understanding the space of Riemannian
metrics on a smooth compact manifold of dimension larger than two; various flows of
metrics have been studied, and several new interesting functionals on this space have been found,
for example via spectral invariants such as functional determinants, or in connection
with conformal geometry: $Q$-curvature, renormalized volumes in AdS/CFT theory, and
fully non-linear equations for certain curvature quantities. As we shall see below, when the
functional is (in addition to being diffeomorphism invariant) conformally invariant,
and the manifold is the $n$-sphere, it becomes remarkably natural to apply representation
theory for the conformal group of the sphere; the (infinite-dimensional) principal
series representations of this group on various tensor fields on the sphere may
be effectively analyzed as Harish-Chandra modules, and in particular one may
calculate explicitly the invariant Hermitian forms on such modules (in general
meromorphic functions of the corresponding parameters, in geometric terms the conformal weights). This takes the form of a spectrum generating principle for determining the invariant object.

Our main results Theorem \ref{universal} and Theorem \ref{universal3} identify (up to a constant) the Hessian of a functional at the standard $n$-sphere (which is a stationary point), with exactly one such Hermitian form (in dimension $n\geq4$, while in dimension 3 the space of such is two-dimensional); thus we may speak of a universal Hessian, and this we identify in even dimensions in Corollary \ref{HEvenDimension} in terms of Fefferman-Graham's obstruction tensor, which is an important object from the theory of conformally compact Poincare-Einstein manifolds.

As our main applications of this rigidity result we obtain (in Theorem \ref{KatesTheorem}), based on knowledge of the leading term in the Hessian, a new and conceptually clear proof of earlier results on extremals of determinants by K. Okikiolu, which appeared first in the paper \cite{Ok1} in Ann. Math. (2001), and also of the natural extensions to general conformally covariant operators (Theorem \ref{confcov}), e.g. to determinants of Dirac operators (Theorem \ref{DiracExtremals}). 

Furthermore we apply our new principle to study the value of the zeta function evaluated at zero in even dimension (Theorems \ref{Hessianzetazero}, \ref{zetazeroDirac} and \ref{zetazeroYamabe}), which has not been studied before. While our approach is not sensitive to the parity of the dimension (as long as the functional is conformally invariant), the even-dimensional case does not seem to be manageable using methods from \cite{Ok1}, since the argument there relied on the odd-dimensional phenomenon of the Kontsevich-Vishik trace. Finally we prove in Theorem \ref{QTheorem} that the total $Q$-curvature has local maxima at the standard spheres.

Note that for instance for the problem of determinants, the conformal covariance of the operator is essential, e.g. it is known that local extremality at the round sphere metrics does not hold for the ordinary Laplacian $\Delta=\nabla^i\nabla_i$ acting on functions (see \cite{Ok1}), namely $\det \Delta$ has a saddle point at $(S^{2k+1},g_0)$ in dimensions $2k+1\geq4$. Some of the virtues of the present paper are to spell out: (1) the role played by the fact that the operator is (an integer power of) a conformally covariant operator, and (2) the importance of the ground metric being the round sphere; our viewpoint in this paper is that $S^n$ has a very large conformal group. With such new insight one realizes that the earlier results by K. Okikiolu for the special case of the determinant of the conformal Laplacian is in fact the quite generic picture for the whole class of conformal functionals, on the $n$-sphere. We think of our results as a step towards a Morse theory for the space of metrics.

The proofs presented here build on Tom Branson's ideas in \cite{BrFuncDet} concerning the role of the complementary series representations of the conformal group of $S^n$, namely $\SO{n+1,1}$, in the search for geometric inequalities relevant for the study of the extremals in a conformal class of the functional determinant in even dimensions (see also \cite{BOO}). Cunningly, as witnessed by the following sections, one is in the conformally invariant case naturally lead to exploiting the action of the conformal group on exactly all of the non-conformal (and non-diffeomorphic) directions. In this connection, Peter Sarnak has noted that the recent paper \cite{SS} (joint with A. Str\"o{}mbergsson) exploits similar ideas in a different setting, where the symmetry group is large but finite.

The spectrum generating principles have been considered previously by many authors (see for instance \cite{BrFuncDet}, and \cite{BOO} for general results), and the present authors learned from R. Graham that an application of the same intertwining operator being identified in this paper appeared in \cite{Gr1}. There it was recently studied as an intertwining operator for different purposes, namely in the context of the Dirichlet-Neumann map for Poincare-Einstein metrics. Also an explicit formula was given there (see Equation (\ref{GrEq})), at least in the odd-dimensional case (and is known in the even-dimensional case, e.g. \cite{Gr2}). It was also shown in the paper \cite{Gr1} that the relevant space of intertwiners on $S^3$ is in fact 2-dimensional in certain cases.

\section{Conformal functionals in Riemannian geometry}
Let $(M^n,g_0)$ be a compact smooth oriented Riemannian manifold (and assume it is spin, with a fixed topological spin structure, whenever needed to define the operators, functionals etc. in our examples). We denote by $\Metr(M)$ the space of all smooth Riemannian metrics on $M$, and study functionals
\[
F:\Metr(M)\to\Reals,
\]
thinking of $g_0$ as the ``ground metric''. Note that $\Metr(M)$ may be realized and given a smooth topology in several ways, e.g. as a Banach or tame Fréchet manifold (see \cite{Eb}, \cite{FrGr}, \cite{Ham}), using Sobolev spaces of sections of the bundle of symmetric two-tensors $S^2TM$, of which the convex cone of Riemannian metrics on $M$ is a subset. The tangent space to the space of smooth Riemannian metrics, at the metric $g_0$, is naturally identified with $\CinfStwoM$ for the details of such constructions). Note that there is a natural $L^2$-integral pairing on two-tensors on $M$, using the metric $g$,
\beq\label{NaturalPairing}
\DIP{h}{k}_g=\int_M\langle h,k\rangle_gdV_g=\int_M h_{ij}k_{lm}g^{il}g^{jm}dV_g.
\eeq
We impose the following geometric assumptions on the functional $F$.
\begin{geoassumptions}\label{geoassump}
$F$ satisfies for any $g\in \Metr(M)$.
\begin{itemize}
\item[(1)] $F(\varphi^*g)=F(g),\quad \varphi\in\DiffO{M}$,\\
\item[(2)] $F(e^{2\omega}g)=F(g),\quad\omega\in\Cinf{M}$.
\end{itemize}
\end{geoassumptions}
\begin{remark}
In dimensions $n\geq4$ it is enough to require, as here, the invariance under the identity component $\DiffO{M}$ of the diffeomorphism group $\Diff{M}$ in order to obtain the rigidity result (for $n=2,3$ see Section \ref{sec:lowdimension}).
\end{remark}
Under these assumptions, and for any $k\in\CinfStwoM$ in the tangent space, $\varphi\in\ctran(M,g_0)$ a conformal transformation of $(M,g_0)$ and $\varphi^*g_0=\Omega^2_\varphi g_0$, with $\Omega_\varphi$ being the corresponding conformal factor, we see that
\beq\label{FundId1}
F(g_0+tk)=F(\varphi^*g_0+t\varphi^*k)=F(\Omega^2_\varphi g_0+t\varphi^*k)=F(g_0+t\Omega^{-2}_\varphi\varphi^*k).
\eeq
Thus the variational problem for the functional $F$ at the ground metric $g_0$ has a certain invariance property under a specific action of the conformal group $\ctran(M,g_0)$ of the manifold $(M,g_0)$.

\begin{analassump}\label{analassump}
$F$ satisfies for each $k\in\CinfStwoM$
\begin{itemize}
\item[(1)] The map $t\to F(g_0+tk)\in\Reals$ is at least $C^3$ in a nbh. of $t=0$, and
\[
D^2F_{g_0}(k,k):=\frac{d^2}{dt^2}_{\big|t=0}F(g_0+tk)=\DIP{k}{H_{g_0}k}_{g_0},
\]
where the Hessian operator $H_{g_0}=H(F,g_0)$ is a linear operator
\[
H_{g_0}:\CinfStwoM\to\CinfStwoM,
\]
which is symmetric with respect to $\DIP{\cdot}{\cdot}_{g_0}$.
\end{itemize}
\end{analassump}
\begin{remark}
In examples it has been verified that the Hessian $H_{g_0}$ exists as an $n$'th order pseudodifferential operator in the vectorbundle $S^2TM$, more specifically a differential operator (see Corollary \ref{HEvenDimension} below), a classical polyhomogeneous pseudodifferential operator (\cite{Ok1}), or a log-polyhomogeneous pseudodifferential operator (see \cite{Mol1}). The pseudodifferential Hessian calculus for zeta functions of geometric Laplace-type operators of order 2 was developed by Kate Okikiolu using heat kernel methods in \cite{Ok1}, \cite{Ok2}, \cite{Ok3}, \cite{OkW}.
\end{remark}

We note that by definition $g_0$ is a stationary point for $F$ if
\[
DF_{g_0}(k):=\frac{d}{dt}_{\big|t=0}F(g_0+tk)=0,\quad\forall k\in\CinfStwoM.
\] 
In fact it follows from \cite{Bl} that any conformal functional must have the standard spheres $(S^n,g_{S^n})$ as stationary points, and thus this will be the case for any of our concrete examples of such functionals.

Exploiting the conformal invariance of $F$, i.e. in the directions tangent to conformal rescalings of $g_0$, given by
\[
\conf_{g_0}:=\big\{\omega g_0\big|\omega\in\Cinf{M}\big\}\subseteq\CinfStwoM,
\]
we shall consider the space $\{\conf_{g_0}\}^\perp=\TFSM$, where the orthogonal complement is with respect to the natural inner product in (\ref{NaturalPairing}). Namely, the full Hessian is given by
\[
H=\begin{pmatrix} \tilde{H} & 0 \\ 0 & 0 \end{pmatrix}:\TFSM\oplus\conf_{g_0}\to\TFSM\oplus\conf_{g_0},
\]
and we will denote the restricted Hessian $\tilde{H}$ simply by $H$.

Note also that from the diffeomorphism invariance, the subspace (where $L_X$ denotes the Lie derivative)
\begin{align*}
\diff_{g_0}:&=\Big\{\fracsm{d}{dt}_{\big\vert t=0}\varphi^*_tg_0\big|\;\textrm{for}\;t\mapsto\varphi_t\in\DiffO{M}\;\textrm{a}\;C^\infty\textrm{-curve s.t.}\;\varphi_0=\Id\Big\}\\
&=\big\{\Lie_Xg_0\big|X\in\Cinf{TM}\big\}\\
&=\big\{k\in\CinfStwoM\big|k(X,Y)=(\nabla_X\omega)(Y)-(\nabla_Y\omega)(X),\;\omega\in\Omega^1(TM)\big\},
\end{align*}
of directions tangent to diffeomorphism pullbacks of $g_0$, is always in the kernel of the Hessian of $F$, and hence so is the restricted space
\[
\diff^0_{g_0}:=\diff_{g_0}\cap\conf_{g_0}^\perp\subseteq\TFSM.
\]

Note that also the sum space $\confdiff$ is in the kernel of the Hessian, for any conformal functional $F$. Furthermore we use
\[
\confdiffperp=\big\{k\in\CinfStwoM\big|\tr_{g_0}k=\diver_{g_0}k=0\big\},
\]
to denote the space of trace- and divergence-free symmetric two-tensor fields, which is the $L^2$-orthogonal complement of $\confdiff$.

Under Assumptions \ref{analassump} we find, upon differentiation of (\ref{FundId1}), that
\beq\label{FormInvariant}
\DIP{\Omega^{-2}_\varphi\varphi^*k}{H_{g_0}\Omega^{-2}_\varphi\varphi^*k}_{g_0}=\DIP{k}{H_{g_0}k}_{g_0},
\eeq
meaning that the Hessian form too is invariant under a certain action, denoted $u_0$, of $\ctran(M,g_0)$. Here $\Omega_\varphi$ appears to the power of $-2$, and we remind that for each $\nu\in\Reals$ there is such an associated group (right-)action. Namely one defines, writing $\rho=\fracsm{n}{2}$
\beq\label{unu}
u_\nu(\varphi)k=\Omega_\varphi^{\rho+\nu-2}\varphi^*k,\quad k\in\TFSM,
\eeq
for $\varphi^*g_0=\Omega^2_\varphi g_0$, where $\varphi\in\ctran(M,g_0)$.
Likewise we define the infinitesimalized version, for $X\in\cvf(M,g_0)$ a conformal vector field and $\omega_X$ the infinitesimal cocycle correspondingly defined by $\Lie_Xg=2\omega_Xg$ as follows:
\[
U_v(X)k=\Lie_Xk+\big(\rho+\nu-2\big)\omega_Xk,\quad k\in\TFSM.
\]
The offsets of ``$-2$'' here are merely conventional, motivated by the fact that the differential geometric realization of $\TFSM$ has internal conformal weight $+2$ (see Section \ref{sec:reptheory}). Likewise the so-called $\rho$-shift appearing above, here $\rho=\fracsm{n}{2}$, will be convenient later.

Since any orientation-preserving $\varphi\in\ctran_0(M,g_0)$ pulls back the Riemannian measure according to
\[
\varphi^*dV=\Omega^n_\varphi dV,
\]
we find by exploiting diffeomorphism invariance of integration on $M$, and $\varphi^*g_0=\Omega_\varphi g_0$ that
\[
\DIP{h}{k}_{g_0}=\DIP{u_{-n/2}(\varphi)h}{u_{n/2}(\varphi)k}_{g_0},
\]
for any $\varphi\in\ctran(M,g_0)$. Using this with (\ref{FormInvariant}) we get
\begin{align*}
\DIP{k}{H_{g_0}k}_{g_0}&=\DIP{u_{-n/2}(\varphi)k}{H_{g_0}u_{-n/2}(\varphi)k}_{g_0}\\
&=\DIP{u_{-n/2}(\varphi^{-1})u_{-n/2}(\varphi)k}{u_{n/2}(\varphi^{-1})H_{g_0}u_{-n/2}(\varphi)k}_{g_0}\\
&=\DIP{k}{u_{n/2}(\varphi^{-1})H_{g_0}u_{-n/2}(\varphi)k}_{g_0},
\end{align*}
for any $k\in\TFSM$, and therefore
\[
H_{g_0}=u_{n/2}(\varphi^{-1})\circ H_{g_0}\circ u_{-n/2}(\varphi).
\]
Thus we have proved the following proposition.
\begin{proposition}
For any conformal functional $F$, the Hessian operator $H(F,g_0)$ satisfies
\beq\label{Hintertwining}
\Omega_\varphi^{n-2}\varphi^*(H_{g_0}k)=H_{g_0}\Omega^{-2}_\varphi\varphi^*k, \quad\text{for}\quad k\in\TFSM.
\eeq
In other words the Hessian of $F$ is an \textit{intertwining operator} between the two representations $u_{-n/2}$ and $u_{n/2}$ of $\ctran(M,g_0)$, and for $U_{-n/2}$ and $U_{n/2}$ of $\cvf(M,g_0)$, on the space of sections of the bundle $\TFSM$ of symmetric trace-free covariant two-tensors.
\end{proposition}

In the light of the proposition, it is quite natural to investigate the representations $u_{-n/2}$ and $u_{n/2}$ appearing here, in order to 
understand the Hessian. These will turn out not to be irreducible, as is to be expected geometrically, since there is an infinite-dimensional kernel of the Hessian (containing at least the natural subspace of $\TFSM$ arising from conformal and gauge invariance), and we know examples of several quantities for which the Hessian is not identically zero (e.g. \cite{Ok1}, \cite{Ok2}, \cite{OkW} and \cite{Mol1}). Studying the quotient space of $\TFSM$ with that natural subspace, one might ask what happens if the conformal group is large, i.e. in the case of the standard spheres, where it is a classical fact that $\ctran(M,g_0)$ has maximal dimension (see e.g. \cite{KN}). And indeed: On $(S^n,g_{S^n})$ an irreducibility result on this quotient or ``moduli space'' does hold. To prove this, we need crucially the geometry of the so-called Ahlfors (or conformal Killing) operator.

\section{Conformal geometry of the Ahlfors operator}\label{sec:conformalgeometry}
\begin{definition}
The Ahlfors operator (or conformal Killing operator)
\[
S_g:\Cinf{TM}\to\TFSM
\]
on $(M,g)$ is defined by
\[
S_gX=\Lie_Xg-\frac{2}{n}(\diver X)g,\quad\in\Cinf{TM}.
\]
\end{definition}

The next proposition lists some basic properties of this important operator. For details and further information, see \cite{OP}.
\begin{proposition}\label{Simker}
\begin{align*}
&\ker S=\cvf(M,g),\\
&S^*=2\diver_g:\TFSM\to\Cinf{TM},\\
&\ran S_g=\big(\big\{k\in\TFSM\big|\diver_g k=0\big\}\big)^\perp=\diff^0_g,
\end{align*}
where $\perp$ is with respect to $\TFSM$.
\end{proposition}
Note that in particular for the Hessian operator, by the discussion in the previous section, we have for any given functional $F$ satisfying our assumptions that
\beq
\ran S_{g_0}\subseteq\ker H(g_0, F).
\eeq
We observe another important property of $S$ in the following proposition.
\begin{proposition}\label{Sconf}
$S$ is conformally covariant, namely for $\varphi\in\ctran(M,g)$,
\beq\label{EqSconf}
\Omega_\varphi^{-2}\varphi^*SX=S\varphi^*X, \quad\textrm{where}\quad\varphi^*g=\Omega^2_\varphi g.
\eeq
\end{proposition}
\Proof
Using that $\Lie_Xg=\frac{d}{dt}_{\big|t=0}\psi_t^*g$, where $\psi_t$ is the flow of $X$,
\beq\label{LieComp}
e^{-2\omega}\varphi^*\Lie_Xg=2(\varphi^*X).\omega g+\Lie_{\varphi^*X}g,
\eeq
where we have written $\omega=\log\Omega_\varphi$. The Levi-Civita connection of the metric pulled back by $\varphi$ is
\[
\nabla_X^{\varphi^*g}Y=\varphi^*\Big(\nabla_{\varphi_*X}^g\varphi_*Y\Big),
\]
and it follows that
\beq
\varphi^*(\diver_g X)=\diver_{\varphi^*g}\varphi^*X=\diver_{e^{2\omega}g}\varphi^*X.
\eeq
As seen from the Koszul formula, $\nabla$ changes conformally according to
\beq\label{NablaConf}
\nabla_X^{e^{2\omega}g}Y=\nabla^g_XY+(X.\omega)Y+(Y.\omega)X-g(X,Y)\nabla\omega.
\eeq
Choosing a local orthonormal frame $\{X_i\}$, we get from (\ref{NablaConf}) that
\beq\label{DivConf}
\diver_{e^{2\omega}}\varphi^*X=\sum_{i=1}^ne^{2\omega}g\bigg(\nabla^{e^{2\omega}g}_{e^{-\omega}X_i}\varphi^*X,e^{-\omega}X_i\bigg)=\diver_g\varphi^*X+n(\varphi^*X).\omega.
\eeq
Combining (\ref{LieComp})-(\ref{DivConf}) the proof is completed:
\[
e^{-2\omega}\varphi^*SX=2(\varphi^*X).\omega g+\Lie_{\varphi^*X}g-e^{-2\omega}\frac{2}{n}\varphi^*(\diver_gX)e^{2\omega}g=S\varphi^*X.
\]
\Endproof

\section{Representation theory of the conformal group, and the universal Hessian}\label{sec:reptheory}
As the manifold and ground metric $(M,g_0)$ we now take the sphere with the standard metric $(S^n,g_{S^n})$, assuming $n\geq2$. We recast the differential geometric notions from the previous sections into the language of representation theory. We refer in particular to \cite{BrFuncDet} and also \cite{BOO} for quite general and detailed expositions of these topics, including many relevant references. The central object is (the identity component of) the conformal group of the round sphere, which is the semisimple Lie group
\[
G=\SOzero{n+1,1}.
\]
The action of $A\in G$ has a geometric realization as follows, by viewing $S^n$ as the unit sphere in $\Reals^{n+1}$ and as $S^n\times\{1\}\subseteq\Reals^{n+2}$. Then for $y\in S^n$,
\[
A\cdot y=\frac{A(y,1)}{(A(y,1))_{n+1}}\in S^n\times\{1\}
\]
defines the action and gives an isomorphism of $G$ with $\ctran(M,g)$.

We fix an Iwasawa decomposition $G=KAN$ relative to a Cartan decomposition $\mathfrak{g}=\mathfrak{k}+\mathfrak{s}$ and maximal abelian Lie subalgebra $\mathfrak{a}\subseteq\mathfrak{s}$, and also fix a minimal parabolic subgroup $P=MAN$. Here $K$ is $\SO{n+1}$ acting in the first $n$ variables, and $M$ is $\SO{n}$ acting in the middle $n$ variables, $A\simeq\Reals$  and $N\simeq\Reals^n$. Recall that
\[
G/P\simeq K/M\simeq S^n.
\]
Half the sum of the positive roots is $\rho=n\alpha/2$, where the $(\mathfrak{g},\mathfrak{a})$ roots are $\pm\alpha$.

If $(V_\sigma,\sigma)\in\Irr(M)$ is an irreducible representation of $\SO{n}$, there is a representation $(a,h)\to a^p\sigma(h)$ for $(a,h)\in\Reals^+\times\SO{n}$. Denoting this by $(V^p_\sigma,\sigma^p)$ one defines associated bundles with conformal weight $p\in\Reals$ by
\begin{align*}
&\mathbb{V}_\sigma^p=G\times_PV_\sigma^p=\Ind_{MAN}^{KAN}\big(\sigma\otimes a^p\otimes 1\big),\\
&\Cinf{\mathbb{V}_\sigma^p}=\bigg\{\psi\in\Cinf{G}\bigg|\psi(xman)=a^{-p-\rho}\sigma(m)^{-1}\psi(x)\bigg\}.
\end{align*}
The differential geometric realization is the bundle associated to the conformal frame bundle $\mathcal{F}_{\Reals^+\times\SO{n}}$ by $\sigma^p$.

The irreducible representations $\sigma\in\Irr(M)$, and analogously for $K=\SO{n+1}$, are parametrized by dominant weight vectors also denoted $\sigma\in\mathbb{Z}^{\lfloor n/2\rfloor}$, with
\beq
\begin{split}
&\sigma_1\geq\ldots\geq\sigma_{[n/2]},\quad n\:\textrm{odd},\\
&\sigma_1\geq\ldots\geq\sigma_{[n/2]-1}\geq\big|\sigma_{[n/2]}\big|,\quad n\:\textrm{even}.
\end{split}
\eeq
Examples are $\sigma=(0)=(0,0,\ldots,0)$ the trivial representations, $\sigma=(1)=(1,0,\ldots,0)$ the defining representation, and $\sigma=(2)=(2,0,\ldots,0)$ the trace free symmetric two tensor representation. The defining representation $(1)$ and conformal weight $+1$ produces the cotangent bundle $T^*M$, while $(1)$ with weight $-1$ gives the tangent bundle $TM$. The realization of $(2)$ as $S^2_0TS^n$ is the bundle with conformal weight $+2$, corresponding to $\mathbb{V}^2_{(2)}$.

The space $\mathcal{E}(S^n,\mathbb{V}^p_\sigma)$ of $K$-finite sections is decomposed into $K$-types using Frobenius reciprocity
\[
\Homo_M(\sigma,\beta_{|M})\simeq\Homo_K(\mathcal{E}(S^n,\mathbb{V}^p_\sigma),\beta),\quad\text{for}\quad\beta\in\Irr(K).
\]
For the multiplicity of $\sigma$ in $\beta_{|M}$ there is a branching law, namely it is either $0$ or $1$, and is $1$ if and only if
\beq\label{branching}
\begin{split}
&\beta_1\geq\sigma_1\geq\beta_2\geq\sigma_2\geq\ldots\geq\sigma_{[n/2]}\geq\big|\beta_{\sigma_{[n/2]+1}}\big|,\quad n\:\textrm{odd},\\
&\beta_1\geq\sigma_1\geq\beta_2\geq\sigma_2\geq\ldots\geq\sigma_{[n/2]-1}\geq\beta_{\sigma_{[n/2]}}\geq\big|\sigma_{[n/2]}\big|,\quad n\:\textrm{even}.
\end{split}
\eeq
When this is the case, we write $\beta\downarrow\sigma$ or $\sigma\uparrow\beta$. Note that the $K$-action is independent of $p$,
\[
\mathbb{V}^p_{\sigma|K}=\Ind_M^K(\sigma).
\]
Using the above we decompose the tangent bundle $\mathbb{V}_{(1)}$ into $K$-types as
\beq\label{TMKtypes}
\mathcal{E}(S^n,\mathbb{V}_{(1)})=\bigoplus_{r=0}^1\bigoplus_{l=0}^{\infty}E_{(1+l,r)},\quad\text{where}\quad E_{(1+l,r)}\simeq_K(1+l,r).
\eeq
Similarly for $\sigma=(2)$ and $n\geq 4$
\beq\label{TFSMKtypes}
\mathcal{E}(S^n,\mathbb{V}_{(2)})=\bigoplus_{q=0}^2\bigoplus_{j=0}^{\infty}F_{(2+j,q)},\quad\text{where}\quad F_{(2+j,q)}\simeq_K(2+j,q), \quad n\geq4.
\eeq
The differential geometric observations from Sections \ref{sec:introduction} and \ref{sec:conformalgeometry} can now be reformulated, where again $\rho=n/2$.
\begin{proposition}\label{reformulation}
Let $F$ be a given conformal functional. The operators $S$ and $H=H(F,g_0)$ are intertwining for the $(\mathfrak{g},K)$-modules
\begin{align*}
&S:\mathbb{V}^{-1}_{(1)}\to\mathbb{V}^0_{(2)},\\
&H:\mathbb{V}^{\rho-\fracsm{n}{2}}_{(2)}\to\mathbb{V}^{\rho+\fracsm{n}{2}}_{(2)}
\end{align*}
In particular $S$ and $H$ are diagonalized by the $K$-decompositions in (\ref{TMKtypes}) and (\ref{TFSMKtypes}), respectively.
\end{proposition}
\Proof
The Ahlfors operator $S$ is intertwining from $\mathbb{V}^{-1}_{(1)}$ to $\mathbb{V}^0_{(2)}$ by Proposition \ref{Sconf}, since the actions on each side in (\ref{EqSconf}) correspond to these conformal weights. The differential geometric realization gives an internal conformal weight of $+2$ on $\TFSM$. Similarly for $H_{g_0}$ by (\ref{Hintertwining}) and the following discussion. That the operators are diagonalized follows from Schur's lemma.
\Endproof
\begin{remark}
We could have given a more direct proof of the intertwining property of $S$, in representation theoretical terms. Namely, $S$ is a generalized gradient and by Fegan's Theorem (see \cite{Fe}, or Theorem 7.5 in \cite{BrFuncDet}), it is intertwining for exactly the conformal weights in Proposition \ref{reformulation}.
\end{remark}
Since $S$ is intertwining, Schur's lemma and the branching laws in (\ref{branching}) gives us the following very crucial result.
\begin{lemma}\label{ranSKtypes}
In terms of $K$-types as above, we have for $n\geq4$,
\[
\ker S=\bigoplus_{r=0}^1E_{(1,r)},\quad\text{and}\quad\ran S=\bigoplus_{q=0}^1\bigoplus_{j=0}^{\infty}F_{(2+j,q)}.
\]
\end{lemma}

Studying now the $(\mathfrak{g},K)$-modules, we turn to the irreducibility issues. The notation is that $\beta\leftrightarrow\gamma$ when $\gamma$ is a $K$-summand of $\mathfrak{s}\otimes\beta$. As the notation suggests, the relation $\beta\leftrightarrow\gamma$ is symmetric (due to $\mathfrak{s}$ being self-dual as a $K$-module). We also write

\begin{align*}
&\kappa_\beta=\left<2\rho_{\text{so}(n)}+\beta,\beta\right>_{\Reals^L},\quad\text{where}\\
&2\rho_{\text{so}(n)}=(n-2,n-4,\ldots,n-2k),\quad\text{if}\quad n=2k\quad\text{or}\quad n=2k+1.
\end{align*}

In particular, if $\beta=(2+j,q)$, for $j\in\mathbb{N}_0$ and $q\in\{0,1,2\}$,
\begin{align}
&k_\beta=(n+j+1)(j+2)+q(n+q-3)\\
&k_{\beta+(1,0)}-k_\beta=n+2j+4,\label{jstepping}\\
&k_{\beta+(0,1)}-k_\beta=n+2q-2.\label{qstepping}
\end{align}

\begin{lemma}[Branson's cocycle irreducibility, \cite{BrFuncDet}, Lemma 7.10]\label{cocycleirred}
Assume the dimension is $n\geq4$. If $\beta\leftrightarrow\gamma$ then
\[
\Proj_\gamma\omega(\mathfrak{s})F_\beta=F_\gamma.
\]
As a result, the span of the orbit of any nonzero $k\in F_\beta$, for any $\beta\downarrow(2)$, under the joint action of $u_\nu(K)$ and $U_\nu(\mathfrak{s})$, is $\mathcal{E}(S^n,\mathbb{V}_{(2)})$. In particular, $\mathcal{E}(S^n,\mathbb{V}_{(2)})$ has no nontrivial invariant subspace under this action.
\end{lemma}

\begin{proposition}[\cite{BrFuncDet}, Corollary 7.11]\label{irred}
Assume the dimension is $n\geq4$. If $\beta=(2+j,q)$ as above, $\beta\leftrightarrow\gamma\downarrow(2)$, and $X\in\mathfrak{g}$,
\beq\label{ProjFormula}
\Proj_\gamma U_\nu(X)_{|\beta}=c(\beta,\gamma,\nu)\Proj_\gamma\omega_{X|\beta},
\eeq
where $c(\beta,\gamma,\nu)=\fracsm{1}{2}(\kappa_\gamma-\kappa_\beta+2\nu)$. If $\big|\nu|<\frac{n}{2}-1$, then
$\Proj_\gamma U_\nu(X)F_\beta=F_\gamma$. As a result, the span of the orbit of any nonzero $\varphi\in F_\beta$, for any $\beta\downarrow(2)$, under the joint action of $u_\nu(K)$ and $U_\nu(\mathfrak{s})$, is $\mathcal{E}(S^n,\mathbb{V}_{(2)})$. In particular, $\mathcal{E}(S^n,\mathbb{V}_{(2)})$ has no nontrivial invariant subspace under this action.
\end{proposition}

Note that because $K$ acts by isometries of the standard $S^n$, the action of $u_\nu(K)$ is in fact $\nu$-independent (while $U_\nu(\mathfrak{g})$ is manifestly not), so we may drop the subscript and write $u(K)=u_\nu(K)$. Since the representations needed in our application have $\nu=\pm\fracsm{n}{2}$, the parameter is outside the domain of irreducibility. As discussed above, this reducibility is a natural consequence of the geometry inherent in the problem. In both of the cases $\nu=\pm\fracsm{n}{2}$ it is seen from (\ref{jstepping}) that stepping in the $j$ direction is always possible, i.e. the map in (\ref{ProjFormula}) is onto. Equation (\ref{qstepping}) shows that in the case of $\nu=-\fracsm{n}{2}$, stepping in $q$ is possible, except for stepping \emph{up from} $(2+j,1)$. Thus irreducibility fails, but by Lemma \ref{ranSKtypes}, Lemma \ref{cocycleirred} and Proposition \ref{irred}, the representation $(u,U_{-n/2})$ does descend to an irreducible representation on the quotient space
\[
\mathbb{V}_{(2)}^0\Big/\ran S.
\]

We move on to apply Proposition \ref{irred} in a spectrum generating argument for the Hessian operator $H$. From the intertwining property of $H$, and its diagonalization on the $K$-types, it follows that
\begin{align*}
\mu_\gamma\:c(\beta,\gamma,-\nu)\Proj_\gamma\omega_{X|\beta}&=H\Proj_\gamma U_{-\nu}(X)_{|\beta}=\Proj_\gamma U_{\nu}(X)H_{|\beta}\\
&=\mu_\beta\:c(\beta,\gamma,\nu)\Proj_\gamma\omega_{X|\beta},
\end{align*}
so that the following relation holds
\beq
\mu_\gamma(\kappa_\gamma-\kappa_\beta-2\nu)=\mu_\beta(\kappa_\gamma-\kappa_\beta+2\nu).
\eeq
Using the step rules in (\ref{jstepping})-(\ref{qstepping}) this determines inductively the eigenvalues, and this is the point of the so-called spectrum generating argument. We have thus finally proved the following.

\begin{theorem}[Rigidity - the Hessian is universal]\label{universal}
Assume $n\geq 4$ and let the functional $F$ satisfy the Assumptions \ref{geoassump} and \ref{analassump}. On the standard spheres $(S^n,g_{S^n})$, the associated Hessian
\[
H(F,g_0):\CinfStwoM\to\CinfStwoM
\]
is
\[
H(F,g_0)=c(F)\cdot T_0,\quad c(F)\in\Reals,
\]
for a constant $c(F)$ depending on the functional $F$, and where $T_0$ is the diagonal intertwining operator given by
\[
T_{0|F_{(2+j,q)}}=\frac{\Gamma(n+j+2)\Gamma(n+q-1)}{\Gamma(j+2)\Gamma(q-1)}\cdot\Id_{F_{(2+j,q)}}\geq0,
\]
on each of the $K$-types $(2+j,q)$.

In particular if $c(F)\neq0$, the Hessian of $F$ is semi-definite, and the kernel is exactly
\[
\ker H=\ran S+\conf_{g_0}=\diff^0_{g_0}+\conf_{g_0}=\diff_{g_0}+\conf_{g_0},
\]
i.e. directions in which $F$ is globally invariant by Assumptions \ref{geoassump}.
\end{theorem}
\begin{remark}
\begin{itemize}
\item[]
\item[(1)]
Note that the eigenvalues of $T_0$ are given here as a meromorphic function in $n$ and $j$, and that if the denominator is $+\infty$, it corresponds to one of the zero eigenvalue spaces for $T_0$.
\item[(2)]
Together with Proposition \ref{Sconf} and Proposition \ref{reformulation}, Theorem \ref{universal} proves that (if the rigidity constant $c(F)\neq0$) there is an exact sequence of conformally covariant operators as follows.
\[
\xymatrix{
\mathbb{V}^{-1}_{(1)}\ar[r]^{S}&\mathbb{V}^0_{(2)}\ar[r]^{H(F)}&\mathbb{V}^{n}_{(2)}.
}
\]
This is a so-called exact BGG-sequence, related to the detour complexes studied by for instance by Branson and Gover (see e.g. \cite{BG}).
\item[(3)]
Operators such as $T_0$ appeared already in the literature (see e.g. \cite{BOO}) and the $T_0$ considered here coincides, up to a constant, with the linearized Dirichlet-to-Neumann operator on $S^n$ considered by Robin Graham (e.g. \cite{Gr1} and \cite{Gr2}). Hence an expression for the operator $T_0$ is known in the following form:

\beq\label{GrEq}
T_0=\begin{cases}
c_0|C|^2, & n=3,\\
c_0W^*W(\Delta+c_1^2)(\Delta+c_2^2)\ldots(\Delta+c_k^2)\sqrt{\Delta+c_{k+1}^2}, & n=2k+5,\\
c_0W^*W(\Delta+c_1^2)(\Delta+c_2^2)\ldots(\Delta+c_k^2), & n=2k+4,
\end{cases}
\eeq
for $k\geq 0$ and $c_i\neq 0$. Here $\Delta=\nabla^*\nabla$ is the rough Laplacian, $W=d\mathscr{W}_{g_0}$ is the linearization at the round metric on $S^n$ of the Weyl curvature viewed as an operator $\mathscr{W}:\Metr(S^n)\to\Cinf{\otimes^4T^*S^n}$, and likewise $C=d\mathscr{C}_{g_0}$ for the Cotton-York tensor $\mathscr{C}$.

\item[(4)]
The formula (\ref{GrEq}) should be compared to the now classical product formulas for the GJMS operators by Tom Branson (in the conformally flat Einstein case, e.g. \cite{Br2}), Robin Graham and Rod Gover (in the general Einstein case, see e.g. \cite{FG2}, \cite{Go1}, \cite{Go2}). In particular it is interesting to ask what role the operator $T_0$, in this more geometric form, might play in the case of a general Einstein manifold.
\end{itemize}
\end{remark}

\section{Rigidity for conformal functionals in the exceptional dimensions 2 and 3}\label{sec:lowdimension}
For dimensions $n = 2$ and $n = 3$ we need to modify the discussion slightly; these cases were not covered
in \cite{BrFuncDet} and the composition series are different - corresponding to the non-uniqueness
of the intertwining operators as observed by R. Graham in the case $n = 3$, see \cite{Gr1}. Here there
arises an additional feature of orientation-sensibility of the functional, which corresponds to the
fact that the rotation group is locally a product: 
$K = \SO{4} \simeq \SO{3} \times \SO{3}$ (local isomorphism).
We first consider this  case, i.e. $n = 3$, which has the relevant set of $K$-types
(same notation as in the previous section)
\[
F_{(2+j,q)}\,, \quad  \text{this time with} \quad j \geq 0, \quad q = 0, \pm 1, \pm 2.
\]
At the conformal weight in question (the parameter $p$ as before), the $K$-types with
$q = 0, \pm 1$ form an invariant subspace, and the quotient is spanned by the remaining
$K$-types with $q = \pm 2$. This quotient is the direct sum of two irreducible summands,
viz. those with $q > 0$ and $q < 0$ respectively. It is precisely this reducibility of the
quotient which explains the two-dimensionality of the intertwining operators, observed by
R. Graham. Here we shall be a little more explicit, in that we give the invariant Hermitian
form, analogous to Theorem 1; the arguments are the same as in the previous section, and
we only have to modify the spectrum-generating argument taking into account the new structure
of the set of $K$-types. The result is as follows.

\begin{theorem}[Rigidity - Universal Hessian theorem in dimension 3]\label{universal3}
 Assume $n = 3$ and let the functional $F$ satisfy the Assumptions \ref{geoassump} and \ref{analassump}. On the standard sphere $(S^3,g_{S^n})$, the associated Hessian
 \[
 H(g_0,F):\CinfStwoM\to\CinfStwoM
 \]
 is
 \[
 H(g_0,F)=c^+(F)\cdot T^+_0 + c^-(F)\cdot T^-_0,\quad c^+(F), c^-(F) \in\Reals,
 \]
 for  constants $c^{\pm}(F)$ depending on the functional $F$, and where $T^{\pm}_0$ are the diagonal intertwining operators given by
 \[
 T^+_{0|F_{(2+j,q)}}=\frac{\Gamma(n+j+2)\Gamma(n+q-1)}{\Gamma(j+2)\Gamma(q-1)}\cdot\Id_{F_{(2+j,q)}}
\quad \text{for}\quad q \geq 0,
 \]
and same formula for $T^-_0$ when $q \leq 0$ on each $K$-type $(2+j,q)$.

 In particular if $c^+(F)c^-(F) < 0$, (the
expression above changes sign between
$q = \pm 2$) the Hessian of $F$ is semi-definite, and the kernel is exactly
 \[
 \ker H=\ran S+\conf_{g_0}=\diff^0_{g_0}+\conf_{g_0}=\diff_{g_0}+\conf_{g_0},
 \]
 i.e. directions in which $F$ is globally invariant by Assumptions \ref{geoassump}.
 \end{theorem}

The condition $c^+(F)c^-(F) < 0$ is verified in concrete cases
 in the same way as for large dimensions, as
explained in the following sections; geometrically it corresponds to the functional
$F$ being insensitive to changes of the orientation.

For $n = 2$, we also need to make only minor changes to the general argument, this
time $K = \SO{3}$, and the relevant $K$-types are
\[
F_{2+j}, \quad\textrm{meaning highest weight}\quad(2+j),\quad\textrm{for}\quad j \geq 0,
\]
occurring with multiplicity two; in this case there is no quotient - corresponding to
the well-known fact that the moduli space (of metrics modulo diffeomorphisms and
conformal changes) consists of a single point in the case of the two-sphere. This
is seen from the representation theory, giving the fact that the Ahlfors operator is
onto, and hence the kernel of the Hessian will be the whole space. Note that the tangent
bundle in this case splits in two (when complexified), and the same for symmetric trace-free
two-tensors; this corresponds to the fact that the action of $\SO{2}$ on $\Reals^2$ complexified is the
sum of two characters. It is hardly surprising that a function on a one-point space
has a universal Hessian, namely zero.               

\section{Local extremals of conformal functionals}
The universality of the Hessian in Theorem \ref{universal} and Theorem \ref{universal3} implies that if the rigidity constant $c(F)>0$ (and analogously when $c(F)<0$) then the Hessian of $F$ at the round sphere metric is positive semidefinite, and positive definite when restricted to the linear subspace $\confdiffperp$.

This immediately gives the following weak extremal result.
\begin{proposition}
Let $F:\Metr(M)\to\Reals$ be a conformal functional on $S^n$, and assume that the corresponding constant $c(F)>0$ is positive. Let $g_t$ be a $C^\infty$-curve of Riemannian metrics such that $g_0=g_{S^n}$ is the round sphere metric, and assume
\[
k:=\frac{d}{dt}_{\big|t=0}g_t\notin\confdiff.
\]
Then there exists $\delta=\delta(F,k)>0$ such that
\beq\label{WeakEx}
0<|t|<\delta\Rightarrow F(g_t)>F(g_0).
\eeq
\end{proposition}
\begin{proof}
Since we evaluate along a smooth curve of metrics, this is Taylor's formula with remainder for $C^3$-functions on a real interval.
\end{proof}

Recall that the Hessian operator $c(F)T_0$ has non-trivial kernel consisting of precisely the gauge directions $\confdiff$. Hence the analysis of local extremals in a neighborhood of $g_0$ could a priori be inconclusive, since the affine space
\[
g_0+\ker T_0=g_0+\confdiff,
\]
does not coincide with the actual, in general curved, set of gauge transformed metrics, along which we have by our basic assumptions that the conformal functional is constant. This issue is resolved by an appropriate change of coordinates, as the next proposition shows.

We have only needed to impose weak analytical assumptions on the functional so far. However, in order to apply an inverse function theorem, one may appropriately realize $\Metr(M)$ as well as the other spaces of sections appearing here, as tame Fréchet manifolds in the sense of Hamilton (see e.g. \cite{Ham}), or as Sobolev spaces where the Banach space version of the inverse function theorem holds (see e.g. Ebin's slice theorem \cite{Eb}). In Hamilton's version of Nash-Moser's theory, the diffeomorphism group $\Diff{M}$ has the structure of a smooth tame Fréchet Lie group. In this setting the following formal treatment can be made rigorous.

\begin{proposition}[Local extremality in $\Metr(M)$]\label{Prop:IFT}
Let
\[F:\Metr(M)\to\Reals
\]
be a conformal functional on $M=S^n$ with sufficient smoothness in the Fréchet (or Sobolev) topology of $\Metr(M)$, and assume that the constant $c(F)>0$ is positive.

Then there exists an open neighborhood $U\ni g_0$ of $g_0=g_{S^n}$ (the round sphere metric), with respect to the Fréchet (or Sobolev) topology on $\Metr(M)$, such that:
\beq\label{LocalEx}
g\in U\subseteq\Metr(M)\Rightarrow F(g)\geq F(g_0),
\eeq
and furthermore
\beq\label{StrictEx}
g\in \big[U\backslash \{g_0\}\big]\cap\big[g_0+\confdiffperp\big]\Rightarrow F(g)>F(g_0).
\eeq
\end{proposition}
\begin{proof}
We first note that the assertion in (\ref{StrictEx}) follows by an application of Taylor's formula for Fréchet spaces, which gives that for $k\in\confdiffperp$
\begin{align*}
F(g_0+k)&=F(g_0)+DF_{g_0}(k)+\frac{1}{2!}D^2F_{g_0}(k,k)+R^{(3)}(k)\\
        &=F(g_0)+\frac{c(F)}{2!}\DIP{T_0k}{k}+R^{(3)}(k),
\end{align*}
since $DF_{g_0}(k)=0$ by criticality of $F$ at $g_0$. Here the remainder term is
\[
R^{(3)}_{g_0}(k)=\frac{1}{2!}\int_0^1(1-t)^{2}D^3F_{g_0+tk}(k,k,k)dt.
\]
Since $T_0$ is positive definite on $\confdiffperp$, the claim follows.

Now, to prove (\ref{LocalEx}) we must change coordinates, and consider the smooth map
\[
\Phi:V_1\times V_2\times V_3\subseteq\Cinf{M}\times\Diffeo(M)\times\confdiffperp\to \Metr(M),
\]
given by
\beq\label{PhiDef}
\Phi(\omega,\phi,k)=e^{2\omega}\phi^*(g_0+k),
\eeq
and where we have restricted to a small enough neighborhood of $(0,\Id,0)$, so that (\ref{PhiDef}) defines a Riemannian metric.

Then the differential mapping into $T_{g_0}\Metr(M)=\CinfStwoM$,
\[
D\Phi(\omega,X,k): \Cinf{M}\times\Cinf{TM}\times\confdiffperp\to \CinfStwoM,
\]
at $(0,\Id,0)$ is given by the expression
\[
D\Phi(\omega,X,k)=2\omega g_0+\Lie_Xg_0+k.
\]
Hence the differential is a bounded linear operator with respect to the Fréchet topology (or in the Sobolev realization, between fixed Sobolev spaces of sections). Furthermore it is bijective between the spaces of smooth sections, since
\[
\confdiff=\ran(2\omega g_0+\Lie_Xg_0).
\]
Invoking the inverse function theorem gives an open neighborhood $V$ of $(0,\Id,0)$ such that $\Phi$ is a diffeomorphism onto $\Phi(V)=U\ni g_0$, and hence for any metric $g\in U$ we have
\[
g=e^{2\omega}\phi^*(g_0+k)\quad\textrm{for some}\quad k\in\confdiffperp,
\]
and hence for any $g\in U$ we get as claimed
\[
F(g)=F(e^{2\omega}\phi^*(g_0+k))=F(\phi^*(k+g_0))=F(g_0+k)\geq F(g_0),
\]
where we have used the invariance properties of the conformal functional, and in the final inequality Equation (\ref{StrictEx}) in the proposition.
\end{proof}

\section{Functional determinants and zeta functions of conformally covariant operators}
The main purpose of this section is the application of the general theory developed in the previous chapters, notably the universal Hessian principle in Theorem \ref{universal}, to obtain analogues of Theorem 2 in \cite{Ok1}. The approach here extends to determinants of a generic (integer) power of a conformally covariant operator. As a concrete example, we apply this to the square of the Atiyah-Singer-Dirac operator $\Dirac^2$.

Note that the extremal problems for determinants considered here are somewhat different from those addressed in \cite{BrFuncDet}. There the results concerned global extremals on $S^6$, but only in conformal directions. From \cite{CY}, \cite{On} and \cite{BrFuncDet} it is true, at least for $n=2,4,6$, that $(-1)^{n/2+1}\det L$ and $(-1)^{n/2}\det\Dirac^2$ are maximized (fixing volume) in the conformal class of the standard sphere exactly when $g$ is a pullback of the standard metric by a conformal diffeomorphism.

Recall the definitions of the spectral zeta function and determinant on a compact manifold for an elliptic partial differential operator $P$ of order $d$ with positive definite (or negative definite) leading symbol, and a real discrete spectrum consisting of eigenvalues $\{\lambda_k\}_{k\in\Naturals}$ with finite multiplicity, and $|\lambda_k|\to\infty$ such that Weyl's law $\lambda_k\simeq k^{-d/n}$ is satisfied. Then the spectral zeta function can by defined as:
\[
\zeta_{|P|}(s):=\sum_{\lambda_k\neq0}|\lambda_k|^{-s},\quad\Real{s}>n/d.
\]
For brevity we shall write $\zeta_P(s)=\zeta_{|P|}(s)$. Under the above assumptions, it follows by using a Mellin transform and the heat kernel expansion, that the zeta function has a meromorphic continuation to $\Complex$ which is regular at zero. In particular one may take the $s$-derivative there, and define the determinant as follows
\[
\det(P):=\exp(-\zeta'_P(0)).
\]

Since from the papers \cite{BO1} and \cite{BO2}, both the determinant and zeta function of a conformally covariant operator (or integer power of such) evaluated at $s=0$ are conformally invariant under suitable assumptions (in odd and even dimension, respectively), we get the following theorems by application of the universal Hessians Theorem \ref{universal} (and (\ref{WeakEx})-(\ref{StrictEx}))

\begin{theorem}\label{confcov}
Let $P$ be an integer power of a conformally covariant operator, with positive leading symbol, $\ker P_{g_0}=0$, and such that with respect to uniform dilations of the metric $P$ has homogeneity degree $-\ord P$. Assume also that $g_0$ is a stationary point of $\det P$, and that the Hessian $H_{g_0}(\det P)$ of the determinant exists (cf. Analytical Assumptions \ref{analassump}) and is not the zero operator.

Then on the odd-dimensional standard spheres $(S^{2k+1},g_{S^n})$
\beq\label{HessDet}
H(\det P)=c(\det P)\cdot T_0,\quad c(\det P)\neq 0.
\eeq
In particular the Hessians are semi-definite, and has precisely
\[
\ker H=\confdiff,
\]
so under the assumptions, $\det P$ assumes either a local maximum or minimum at $(S^{2k+1},g_{S^n})$.
\end{theorem}

\begin{theorem}\label{Hessianzetazero}
Under the analogous assumptions on $P$ for the functional $F(g)=\zeta_{P_g}(0)$, we have on the even-dimensional standard spheres $(S^{2k},g_{S^n})$ that
\beq
H(\zeta_P(0))=c(\zeta_P(0))\cdot T_0,\quad c(\zeta_P(0))\neq 0.
\eeq
In particular these Hessians are semi-definite, and has precisely
\[
\ker H=\confdiff,
\]
so under the assumptions $\zeta_P(0)$ assumes either a local maximum or minimum at $(S^{2k},g_{S^n})$.
\end{theorem}

\begin{remark}
\begin{itemize}
\item[]
\item[(1)] Note that all extremals here are strict, apart from in the directions $\confdiff$ (corresponding to the globally invariant directions).
\item[(2)] While the object $T_0$ is always defined, the relations in (\ref{HessDet}) cannot be extended to all dimensions, since the determinant (respectively $\zeta_P(0)$) is not always conformally invariant. This reflects the subtle relation of conformal invariance with the parity of the dimension. Furthermore the Hessian of the determinant is a log-polyhomogeneous operator in some even-dimensional examples (see \cite{Mol1}). The situation for other functionals is however different (see Theorem \ref{QTheorem} below).

\end{itemize}
\end{remark}
Thus, under the assumptions in Theorem \ref{confcov} (that imply Assumption \ref{geoassump}), the determinant $\det P_g$ has local extremals of a type determined by the rigidity constant $c(F)$ in Theorem \ref{universal}, which in turn depends on the operator $P$ in question. To obtain the constant $c(F)$ needed for determining whether the extremum is a maximum or minimum may still constitute a substantial amount of work. In the papers \cite{Ok1}, \cite{Ok2}, \cite{Ok3} and \cite{OkW} the framework for finding the leading symbol of the Hessian operator for zeta functions and determinants of Laplace-type operators of order 2 has been developed. The method is to express the zeta function as a Mellin transform and use heat kernel expansions, the main term under investigation being
\[
\int\int_{u+v<1}(u+v)^sP_k'e^{-uP}P_k'e^{-vP}dudv,
\]
where $P_k'$ is the derivative at $t=0$ of the operator $P$ along a curve of metrics $g+tk$. The analysis in \cite{Ok2} yields the combinatorially quite complicated explicit formula for the leading symbol of the Hessian, in normal geodesic coordinates around a point, in terms of the expression of $P'$ in local coordinates, in the case $\ord P=2$.

For the determinant $\det L$ of the conformal Laplacian,
\[
L:=-\Delta+\frac{n-2}{4(n-1)}\Scal,
\]
where $\Scal$ is the scalar curvature and $\Delta=\nabla^i\nabla_i$ is the ordinary (connection) Laplacian, the explicit leading symbols appear in \cite{Ok1}. Alternatively one may compute for the zeta-function itself, using \cite{Ok2} and variation formulas for the scalar curvature (which has also been done in \cite{Ok3}) to find that the leading symbol of the Hessian of $\zeta_L(s)$ is given by
\begin{align*}
&\langle k,\sigma_{n-2s}(x,\xi)\,k\rangle_g=\Big(\frac{1}{4\pi}\Big)^{\fracsm{n}{2}}\frac{\Gamma(s-n/2)\Gamma(-s+n/2+1)^2}{\Gamma(s)\Gamma(-2s+n+2)}\\
&\times|\xi|^{n-2s}\Bigg\{\Big(\frac{s^2}{(n-1)^2}-\frac{s}{(n-1)^2}-\frac{1}{2}\frac{1}{n-1}\Big)\big(\tr{K_g\Pi^{\bot}_\xi\big)^2}+\frac{1}{2}\tr\big(K_g\Pi^{\bot}_\xi\big)^2\Bigg\},
\end{align*}
for $\Real{s}<n/2-1$, where $\Pi^{\bot}_\xi$ is the orthogonal projection on ${\xi}^{\perp}$, for $\xi\in T_x^{*}M$.

Therefore the sign of $c(\det L)$ is $(-1)^k$, and we obtain a new proof of the following theorem.

\begin{theorem}[\cite{Ok1}]\label{KatesTheorem}
Among metrics on $S^n$ of fixed volume, the standard sphere $(S^{2k+1},g_{S^n})$ is a local maximum for $(-1)^{k+1}\det L$.
\end{theorem}

Furthermore we obtain the following new theorem:

\begin{theorem}\label{zetazeroYamabe}
Among metrics on $S^n$ of fixed volume, the standard sphere $(S^{2k},g_{S^n})$ is a local maximum for $(-1)^{k+1}\zeta_L(0)$.
\end{theorem}

In \cite{Mol1} the Hessian calculus is extended to the case of the square of the Atiyah-Singer-Dirac operator, for general variations of the metric (with a fixed spin structure). The main theorem there is the following, giving the leading symbol of the Hessian for the zeta-function in the meromorphic parameter $s$.

\begin{theorem}[\cite{Mol1}]\label{NielsDirac}
Let $(M^n,\gamma)$ be a closed Riemannian spin manifold. Assume that the kernel of its
Atiyah-Singer-Dirac operator $\Dirac$ has stable
dimension under local variations of the metric, with fixed topological
spin structure.
\newline\indent{}Then the Hessian of the zeta function $\zeta(s)$ of $\Dirac^2$ is a pseudodifferential operator, with leading symbol given by
\begin{align*}
\langle k,\sigma_{n-2s}(x,\xi)\,k\rangle_g&=2^{\lfloor\frac{n}{2}\rfloor-2}\Big(\frac{1}{4\pi}\Big)^{\fracsm{n}{2}}\frac{\Gamma(s-n/2)\Gamma(-s+n/2+1)^2}{\Gamma(s)\Gamma(-2s+n+2)}\\
&\times|\xi|^{n-2s}\Bigg\{\Big[2s-(n-1)\Big]\tr\big(K_g\Pi^{\bot}_\xi\big)^2
+\big(\tr{K_g\Pi^{\bot}_\xi\big)^2}\Bigg\}
\end{align*}
for $\Real{s}<n/2-1$, where $K_g$ is the endomorphism associated to $k$ by raising an index with $g$, and $\Pi^{\bot}_\xi$ is the orthogonal projection on ${\xi}^{\perp}$, for $\xi\in T_x^{*}M$.
\end{theorem}

From this we obtain, by differentiation in $s$, that the sign of the relevant constant $c(\det\Dirac^2)$ is $(-1)^{k+1}$. Furthermore the standard spheres are stationary points of the zeta function $\zeta_{\Dirac^2}(s)$ at each point (see \cite{Mol1} for details), and hence also of $\det P$. Thus we have proved the following theorem.

\begin{theorem}\label{DiracExtremals}
Among metrics on $S^{2k+1}$ of fixed volume, the standard sphere $(S^{2k+1},g_{S^n})$ is a local maximum for $(-1)^k\det\Dirac^2$.
\end{theorem}
Also from Theorem \ref{NielsDirac} we obtain that the sign of the constant $c(\zeta_{\Dirac^2}(0))$ is $(-1)^{k+1}$, thus
\begin{theorem}\label{zetazeroDirac}
Among metrics on $S^{2k}$ of fixed volume, the standard sphere $(S^{2k},g_{S^n})$ is a local maximum for $(-1)^k\zeta_{\Dirac^2}(0)$.
\end{theorem}

\begin{remark}
For yet more examples of such alternating behavior modulo 4 in the dimension of the manifold, for zeta regularized quantities, we point the reader to \cite{Mol2} which deals with the explicit values of the determinant at the stationary points discussed above. E.g. the sign of $\log\det(\Dirac^2,S^n)$ is $(-1)^{\lfloor(n-1)/2\rfloor}$, and $\lim_{n\to\infty}\det(\Dirac^2,S^n)=1$.
\end{remark}

\section{Criteria for the types of extremals of determinants}
It is worth exploring alternative ways of determining the sign of the constant $c(\det P)$, which do not rely on the above quite involved leading symbol calculus, which indeed extracts much more information than needed for our purpose. For instance for $(S^3,g_{S^n})$ the sign of $c(\det\Dirac^2)$ can be found (\cite{Mol3}), by evaluating the second derivative of $\det\Dirac^2$ along a family of Berger metrics on $S^3=SU(2)$, using \cite{Hi}. Since such deformations are non-diffeomorphic and non-conformal and the second variation has a definite sign, this gives together with Theorem \ref{universal} a different proof of the $n=3$ case of Theorem \ref{DiracExtremals} above. While of course $S^{2k+1}$ does not in general have a Lie group structure, it is likely that a similar approach should work for higher dimensions using instead C. B\"ar's construction of Berger-type metrics in \cite{Ba} and the explicit formulas for the eigenvalues of the operator $\Dirac$ derived there.

Furthermore following \cite{Ok1} we remind that one may use the Kontsevich-Vishik trace $\TR$ to express the Hessian form.  We shall not recall all details of this construction here, but merely comment that both the conformal Laplacian (i.e. Yamabe operator) and the Dirac operator have the property that the Green's function $P^{-1}$ is purely singular, so that $\TR DP^{-1}=0$ for any differential operator $D$. Given this, one may write
\[
\DIP{k}{H(\log\det P)k}=-\TR P_k'P^{-1}P_k'P^{-1},
\]
where $P_k'$ is the first variation of $P_g$ in the direction $k\in\CinfStwoM$. For the Dirac operator this (non-conformal) metric variation is to be understood properly using the Bourguignon-Gauduchon formulas, as discussed in \cite{Mol1}. We then obtain the following criterion, emphasizing the fact that one needs very little information about the variational problem to determine the extremum types.
\begin{corollary}
If there exists $k_0\in\CinfStwoM$ such that
\[
\eta(k_0):=-\TR P_{k_0}'P^{-1}P_{k_0}'P^{-1}\neq0,
\]
then the extrema in Theorem \ref{confcov} are maxima when $\eta(k_0)<0$, and minima when $\eta(k_0)>0$.
\end{corollary}
An interesting observation concerning the trace in this corollary is that the much simpler expression $\TR P^{-2}$ in our main examples of the Dirac operator $\Dirac$ and conformal Laplacian $L$, turns out to have the appropriate sign in any dimension.

To evaluate $\TR L^{-2}$, one must derive an expression for the Green's function $G_{L^2}$ of the fourth order operator $L^2$, i.e. the integral kernel of $L^{-2}$. In the following we shall exploit the rotational invariance of $L$ and hence of the Green's function. Fixing a point $y$ on $S^n$ write the conformal Laplacian $L$ in polar geodesic coordinates $(r,\theta)\in[0,\pi)\times S^{n-1}$ around $y$ as
\beq
Lf=-\partial^2_rf-(n-1)\frac{\cos r}{\sin r}\partial_rf+\frac{n(n-2)}{4}f,
\eeq
acting on radial functions $f=f(r)$. Furthermore we see that writing $G_L(x,y)=G_L(r)$ the following equation must hold
\beq\label{Ltwoeq}
LG_{L^2}(r)=G_L(r),\quad\textrm{for}\quad r\neq 0.
\eeq
For the equation $LG_{L}(r)=\delta_0$, and likewise (\ref{Ltwoeq}) the relevant second order linear ODE to study is
\beq
-\partial^2_rf-(n-1)\frac{\cos r}{\sin r}\partial_rf+\frac{n(n-2)}{4}f=0.
\eeq
Using the substitution $z=\cos r$ transforms this to into the form
\beq
(1-z^2)y''-nzy'-\frac{n(n-2)}{4}y=0,
\eeq
which is a hypergeometric equation with the full solution
\beq\label{FullHom}
y=A(1-z)^{-\frac{n-2}{2}}+B(1+z)^{-\frac{n-2}{2}}.
\eeq

From the requirement that the Green's function must be regular at $r=\pi$, that is at $z=-1$, it follows that $B=0$, and by normalization one finds that
\beq
G_{L}(r)=\frac{C_n}{\sin^{n-2}(\frac{r}{2})},\quad\textrm{where}\quad C_n=\frac{1}{2^{n-1}(n-2)\omega_{n-1}}.
\eeq
To determine $G_{L^2}$ we rewrite equation (\ref{Ltwoeq}) as
\beq\label{InHomEq}
(1-z^2)y''-nzy'-\frac{n(n-2)}{4}y=-D_n(1-z)^{-\frac{n-2}{2}},
\eeq
where $D_n=\frac{1}{2^{n/2}(n-2)\omega_{n-1}}$. Using (\ref{FullHom}) and a computation of the Wronskian $W=(1-z^2)^{-n/2}$ leads to the full solution of (\ref{InHomEq})
\beq
y=A(1-z)^{-\frac{n-2}{2}}+B(1+z)^{-\frac{n-2}{2}}-D_n(1+z)^{-\frac{n-2}{2}}\int\Big(\frac{1-z}{1+z}\Big)^{-\frac{n-2}{2}}dz.
\eeq
Requiring again regularity at $z=-1$ we see that with $r=\cos z$
\beq
G_{L^2}(r)=A(1-z)^{-\frac{n-2}{2}}-D_n(1+z)^{-\frac{n-2}{2}}\int_{-1}^z\Big(\frac{1-w}{1+w}\Big)^{-\frac{n-2}{2}}dw.
\eeq
Changing variables to $\tau=\frac{1-w}{1+w}$ and $|x|^2=\frac{1-\cos r}{1+\cos r}$ allows us to express the integral here as follows for $n=2k+1$
\begin{align*}
&\int_{-1}^z\Big(\frac{1-w}{1+w}\Big)^{-\frac{n-2}{2}}dw\\
&\quad=\frac{1}{4}\int_{|x|}^{\infty}\frac{1}{(1+\tau^2)^2}\frac{1}{\tau^{n-1}}d\tau\\
&\quad=\frac{\tau^{2-n}}{4}\prefix_2F_1(1-\fracsm{n}{2},2;2-\fracsm{n}{2};-\tau^2)\Big|_{\tau=|x|}^{\tau=\infty}\\
&\quad=\frac{(-1)^{k}}{4}\Big[\frac{n\pi}{4}-\frac{n}{2}\arctan(|x|)-\frac{|x|}{(1+|x|^2)^2}-\sum_{j=0}^{k-1}\frac{(-1)^j(k-j)}{2j+1}|x|^{-2j-1}\Big].
\end{align*}
From this expansion we can extract the regular part $r\to0$ of the Green's function (in the sense of \cite{Ok1}), which is
\[
G_{L^2}^\textrm{reg}(x,x)=\frac{(-1)^{k+1}(2k+1)\pi}{2^{2k+4}(2k-1)\omega_{2k}},
\]
and thus finally on $(S^n,g_0)$ the value of the Kontsevich-Vishik trace is
\[
\TR L^{-2}=\int_{S^n}G^\textrm{reg}_{L^2}(x,x)dVol(x)=(-1)^{k+1}\frac{\pi^2}{2^{4k+4}}\frac{(2k+1)(2k)!}{(2k-1)(k!)^2}.
\]

For $\Dirac^2$ on $S^n$ the Green's function can be found for instance in the reference \cite{Fi}, and the integral in the representation there can easily be performed similarly to the above, for $n=2k+1$
\begin{align*}
G_{\Dirac^2}(x,y)&=\frac{1}{\omega_{n-1}}\Big(\frac{4}{1+|x|^2}\Big)^{\frac{1-n}{2}}\int_{|x|}^\infty\frac{2}{1+\tau^2}\frac{1}{\tau^{n-1}}d\tau,\\
&=\frac{(-1)^{k}}{\omega_{n-1}}\Big[\frac{\pi}{2}-\arctan(|x|)-\sum_{j=0}^{k-1}\frac{(-1)^j}{2j+1}|x|^{-2j-1}\Big],
\end{align*}
where $1/|x|$ is the the radial coordinate in the stereographic projection from $y$, and $\omega_{n-1}$ is the volume of $S^{n-1}$. Extracting again the regular part and integrating over $S^n$ one obtains
\[
\TR\Dirac^{-2}=(-1)^k\frac{\pi^2}{2^{2k+1}}\frac{(2k!}{(k!)^2}.
\]
We have thus proved the following proposition.
\begin{proposition}
In odd dimension $n=2k+1$ we have, for the Dirac operator and conformal Laplacian, respectively, on the round spheres $(S^n,g_{S^n})$ that
\begin{align*}
&\sign(\TR L^{-2})=(-1)^{k+1},\\
&\sign(\TR\Dirac^{-2})=(-1)^k.
\end{align*}
\end{proposition}

On the basis of this proposition, we offer the following conjecture.
\begin{conjecture}
If $n=2k+1$ and $P$ is a conformally covariant operator, then on the round sphere $(S^n,g_{S^n})$
\[\label{TRConjecture}
\sign(\TR P^{-2})=\sign(\TR P_{k_0}'P^{-1}P_{k_0}'P^{-1})\neq0,
\]
for some and hence all $k_0\in\TFSSn\backslash\confdiff$.
\end{conjecture}
Note that Conjecture \ref{TRConjecture} expresses in particular that all information about the types of the local extremals on $S^n$, near the round metric, should be contained in just the Green's function at the round sphere metric itself. With a criterion such as this, one may readily apply our rigidity theorem to prove extremality of $\det P$ for (integer powers of) any given conformally covariant operator of order $m$, as for instance the GJMS operators \cite{GJMS}. In principle $\TR P^{-2}$ can be computed on $S^n$ as above, by solving a $m$'th order linear ODE involving the explicit Green's function of $P$ for the round sphere metric, and extracting the regular part.

\section{Local extremals of the total $Q$-curvature}
An important theorem in conformal geometry states that the total $Q$-curvature of the Riemannian manifold $(M^n,g)$ for $n=2k$ even,
\[
L(g):=\int_M Q_g dv_g,
\]
is conformally invariant (cf. Theorem 3 in \cite{GZ}, and also \cite{GH}). Since it is also a natural Riemannian invariant, it satisfies our Assumptions \ref{geoassump}, and we may apply our rigidity principle to this functional. In this section we prove the following theorem.

\begin{theorem}\label{QTheorem} Let $n=2k\geq 4$. The total $Q$-curvature has a local maximum at the round sphere metric $(S^{n},g_{S^n})$, under general variations of the Riemannian metric. Namely, there exists an open neighborhood $U\subseteq\Metr(M)$ of $g_{S^n}$ (as in (\ref{LocalEx})) such that
\[
\int_{S^n} Q_{g} dv_g\leq\int_{S^n} Q_{g_{S^n}} dv_{g_{S^n}},\quad\textrm{for any}\quad g\in U\subseteq\Metr(M).
\]
\end{theorem}
\begin{remark}
It follows furthermore from the proof below, that for any closed, conformally Einstein manifold $(M^{2k},g_0)$ there exists a finite-dimensional subspace
\[
V_{(M,g_0)}\subseteq\CinfStwoM,
\]
with possibly $V=\{0\}$ as in the round sphere case, such that the total $Q$-curvature has a strict local maximum at the metric $g_0$, apart from in the directions $k\in V_{(M,g_0)}+\confdiff$ (as in (\ref{WeakEx})).
\end{remark}

See the reference \cite{Han} for related results on the non-critical $Q$-curvatures $Q_{2m}$. There it was proved that the corresponding (rescaled) total $Q$-curvatures are maximized in the conformal class of the standard metric on $S^n$, when the dimension $n$ is odd and $m=\frac{n+1}{2}$ or $\frac{n+3}{2}$. Moreover for $n$ odd and $m\geq\frac{n+5}{2}$ the standard metric on $S^n$ is not stable.

From the proof of our theorem below, it will be clear that Assumptions \ref{analassump} are satisfied. The starting point for the study of extremals of the total $Q$-curvature is the variational formula by Graham and Hirachi. In the following $P$ will denote the Schouten tensor
\[
P=\frac{1}{n-2}\bigg\{\Ric-\frac{\Scal}{2(n-1)}g\bigg\},
\]
where $\Ric$ and $\Scal$ are the Ricci tensor and the scalar curvature in the metric $g$.
\begin{theorem}[Graham-Hirachi's first variation formula \cite{GH}]
Let $M$ be a compact manifold of even dimension $n\geq4$, then for a smooth family of metrics $g_t$,
\beq\label{GHFirstVariation}
\bigg(\int_MQdv\bigg)^{\bullet}=(-1)^{n/2}\frac{n-2}{2}\int_M\langle\mathcal{O}_g,\dot{g}\rangle_g,
\eeq
where the two-tensor $\mathcal{O}$ is the Fefferman-Graham obstruction tensor.
\end{theorem}
We need the detailed properties of the Fefferman-Graham obstruction tensor $\mathcal{O}$.
\begin{theorem}[Properties of Fefferman-Graham's tensor \cite{GH}]\label{OProperties}
The obstruction tensor $\mathcal{O}$ of the even-dimensional Riemannian manifold $(M^n,g)$ has the following properties:
\begin{itemize}
\item[(1)] $\mathcal{O}$ is a natural tensor invariant, i.e. the components are in local coordinates given by universal polynomials in the components of $g$, $g^{-1}$ and the curvature tensor of $g$ and its covariant derivatives.
\item[(2)] $\mathcal{O}$ is symmetric, trace- and divergence-free, and of the form
\beq\label{Oformula}
\mathcal{O}_{ij}=\Delta^{n/2-2}\big(P_{ij,k}{}^k-P_k{}^k,_{ij}\big)+\textrm{LOTS},
\eeq
where $\Delta=\nabla^i\nabla_i$ is the connection Laplacian, and LOTS denotes quadratic and higher terms in curvature involving fewer derivatives.
\item[(3)] If $g$ is conformally Einstein, then $\mathcal{O}_{ij}=0$.
\end{itemize}
\end{theorem}

We can now explain the proof of Theorem \ref{QTheorem}.
\Proof[Proof of Theorem \ref{QTheorem}]
From formula (\ref{GHFirstVariation}) and property (3) in Theorem \ref{OProperties}, we see in particular that the standard metric on the sphere $(S^{2k},g_{S^n})$ is always a stationary point, under arbitrary variations of the metric.

As we will see below, Theorem \ref{OProperties} will ensure that the Hessian of the total $Q$-curvature $L(g)$ exists as a differential operator
\[
H:\CinfStwoM\to\CinfStwoM.
\]
Applying the universal Hessian theorem for $(S^{n},g_{S^n})$ in Theorem \ref{universal} there exist constants $c_n(L)$ such that for each even dimension $n=2k$,
\[
H(L)=c_n(L)\cdot T_0,
\]
where $T_0$ is positive semi-definite and with $\ker T_0$ equal to precisely the trace- and divergence-free symmetric two-tensor fields, corresponding to the global invariance directions of the total $Q$-curvature $L$. According to whether $c_n(L)>0$, $c_n(L)<0$ or $c_n(L)=0$ the analysis shows respectively a local minimum, a local maximum or is inconclusive. In the two first cases, the extremals are furthermore strict, apart from in the space of gauge and conformal invariance directions which equals $\ker T_0$ (see (\ref{WeakEx}) and Proposition \ref{Prop:IFT}).

To find the signs of the constants $c_n(L)$, we compute the leading order symbol of $H(L)$, which involves studying the first variation $\dot{\mathcal{O}}$ of the obstruction tensor.

For a 1-parameter family of metrics $g_t$ and correspondingly natural tensor fields $A_t$ in the sense discussed above, and any operator $\mathcal{C}$ that groups indices of $A_t$ into pairs and performs contractions using the metric $g_t$, we have
\begin{align}
&\big(\nabla^lA\big)^{\bullet}=\nabla^l\dot{A}+\textrm{LOTS},\label{NablaDot}\\
&\big(\mathcal{C}(A)\big)^{\bullet}=C(\dot{A})+\textrm{LOTS},\label{ContractionDot}
\end{align}
where LOTS are terms with a lower total number of spatial and $t$-derivatives in the components of $A$. This notion of lower order terms is consistent with the one used above in (\ref{Oformula}), and is preserved under composition with differential operators, including differentition in $t$ along $g_t$. The point of this leading order calculus is that we only need to determine (the sign of) the leading symbol of the Hessian as a differential operator acting on the section $k\in\CinfStwoM$, and this is given by the terms with the highest number of spatial derivatives of $k$.

In general, if $U_g$ is a local scalar Riemannian invariant, we see that for $g_t=g+tk$, with $k$ trace-free, the first variation at $t=0$ of the integrated invariant is
\[
\bigg(\int_MUdv\bigg)^{\bullet}=\int_M\big(\dot{U}_g+\fracsm{1}{2}U_g\tr_g\dot{g}\big)dv_g=\int_M \dot{U}_gdv_g,
\]
since $\tr_g\dot{g}=0$. Note also that by (\ref{ContractionDot}) we have $(\langle T,\dot{g}\rangle)^{\cdot}=\langle \dot{T}_g,\dot{g}\rangle_g+\textrm{LOTS}$ for any natural 2-tensor. Thus for any trace-free $k\in\CinfStwoM$ and $g_t=g+tk$,
\beq\label{IntQDot}
\bigg(\int_MQdv\bigg)^{\bullet\bullet}=(-1)^{n/2}\frac{n-2}{2}\int_M\big(\langle\dot{\mathcal{O}}_g,k\rangle_g+\textrm{LOTS}\big).
\eeq
Using the naturality in (1) from Theorem \ref{OProperties} and partial integrations, we see that the Hessian here exists as an order $n$ partial differential operator $H$, such that
\[
\frac{d^2}{dt^2}_{|t=0}L(g+tk)=\int_M{\langle Hk,k\rangle}_gdv_g,\quad\textrm{for}\quad k\in\CinfStwoM.
\]

Note that rewriting (\ref{Oformula}) and taking the $t$-derivative using (\ref{NablaDot}) and (\ref{ContractionDot}) gives
\beq\label{ODot}
\dot{\mathcal{O}}=\Delta^{n/2-2}\bigg(\Delta\dot{P}-\frac{1}{2(n-1)}\nabla^2\dot{\Scal}\bigg)+\textrm{LOTS},
\eeq
where $\nabla$ and $\Delta$ are again the connection and connection Laplacian.
The first variation of the scalar curvature is as follows (see e.g. \cite{CLN} or \cite{MT} for a good reference for such computations).
\beq\label{ScalDot}
\frac{\partial\Scal}{\partial t}_{|t=0}=-\langle\Ric,k\rangle_g+\diver_g(\diver_g k)-\Delta\tr_gk,
\eeq
For the Ricci tensor one has
\beq
\frac{\partial\Ric}{\partial t}_{|t=0}=-\frac{1}{2}\Delta_Lk-\frac{1}{2}\Lie_{(\diver_g k)^{\#}}g-\frac{1}{2}\Hess(\tr_gk).
\eeq
Here $\Delta_L$ is the Lichnerowicz Laplacian
\[
(\Delta_Lk)(X,W)=(\Delta k)(X,W)+2\tr k(R(X,\cdot)\cdot,W)-k(X,\Ric(W))-k(W,\Ric(X)),
\]
where $\Delta k$ denotes the connection Laplacian of the two-tensor $k$, $R$ is Riemann's curvature tensor and $\Ric(X)=(\Ric(X,\cdot)^{\#}$.

Restricting to trace- and divergence-free fields $k\in\CinfStwoM$, we obtain
\beq\label{SchoutenVar}
\dot{P}=\frac{1}{2(n-2)}\bigg\{-\Delta_Lk+\frac{1}{n-1}\langle\Ric,k\rangle_gg\bigg\}=-\frac{1}{2(n-2)}\Delta k+\textrm{LOTS}.
\eeq
By combining Equations (\ref{ODot}), (\ref{ScalDot}) and (\ref{SchoutenVar}) we arrive at the result
\[
\dot{\mathcal{O}}=-\frac{1}{2(n-2)}\Delta^{n/2}k+\textrm{LOTS}.
\]
Finally, by (\ref{IntQDot}) the leading symbol of the $n$'th order Hessian operator $H$ is given by
\beq\label{TotalQSymbol}
\sigma_n(H)(x,\xi)k=-\frac{(-1)^{n/2}}{4}\sigma_n\big(\Delta^{n/2}\big)k=-\frac{|\xi|^n}{4}k\Id\in\Endo{S^2TM_x},
\eeq
which holds whenever $k$ is trace and divergence free.

In the presence of the gauge invariance, namely on the space of fields $k\in\confdiff$ that are orthogonal to the trace- and divergence-free fields, ellipticity of the symbol does not hold in the corresponding tangent directions. Therefore one needs to use for instance a factorization of the Hessian as constructed in \cite{Mol1}. Thus
\[
H(L)=\Pi_{\confdiffperp}\tilde{H}\Pi_{\confdiffperp},
\]
where $\Pi$ denotes the $L^2$-orthogonal projection onto the subspace $\confdiffperp$, and $\tilde{H}$ is now a classical pseudodifferential operator, which is elliptic with negative definite leading symbol, viz. the symbol computed above in (\ref{TotalQSymbol}). Then by standard elliptic theory on a closed manifold, applied to the new operator $\tilde{H}$, $H$ is upper semi-bounded, ensuring that $c_n(L)<0$ for the rigidity constant in the round sphere case.

Furthermore, as remarked after the statement of the theorem, this implies that on any $(M^n,g_0)$ which is conformally Einstein (and thus as remarked a stationary point under all variations), the total $Q$-curvature has a local maximum at $g_0$, apart from possibly in a finite number of directions that define a subspace $V_{(M,g_0)}$, namely the direct sum of the finitely many finite-dimensional eigenspaces for the non-negative eigenvalues of $\tilde{H}$.
\Endproof

As an interesting insight from the study of the concrete example of total $Q$-curvature, which for our purposes constitutes a particularly well-understood example of a conformal functional, we have the following corollary, which follows from the proof of Theorem \ref{QTheorem}.
\begin{corollary}\label{HEvenDimension}
In even dimensions, the universal Hessian operator $T_0$ is a natural differential operator, which is given up to a constant as a Hermitian form by the second variation of the total $Q$-curvature, or equivalently through the first variation of the obstruction tensor.
\end{corollary}
Of course this is also clear from Graham's explicit formula in Equation (\ref{GrEq}). In dimensions 4 and 6 where explicit formulas for the $Q$-curvature are known, it is possible to calculate from the above the explicit expression for the universal Hessian operator. Furthermore from the recent developments around the holographic formula for $Q$-curvature found by Andreas Juhl and Robin Graham (see \cite{JG} and \cite{Ju}), one may in principle obtain a scheme for computing the operator $T_0$ in terms of the curvatures and its covariant derivatives. Again, connecting this to Graham's formula (\ref{GrEq}) makes this observation somewhat redundant.

\section{Further examples of conformal functionals}
There are several other natural examples of functionals that fit our theorem. Notably the total $Q$-curvature studied in the previous section furnishes just a specific example of functionals of the type
\[
F(g):=\int_MA(g)dv_g,
\]
where $A(g)$ is a local Riemannian scalar invariant, with the property that $F$ is conformally invariant. Exactly such functionals are the topic of Deser-Schwimmer's Conjecture, which is a statement about the structure of the $A(g)$ (under the conformal invariance assumption on $F$) proved recently by Spyros Alexakis (\cite{Al1}, \cite{Al2}, \cite{Al3}).

Yet another type of conformal functional studied in the literature, comes about by taking the infimum or supremum over conformal classes. One such example is the functional
\[
\tau(g)=\inf_{g\in{[g_0]}}\lambda_1(L_g)\Vol(M,g)^{2/n},
\]
where $\lambda_1(L_g)$ denotes the first eigenvalue of the conformal Laplacian $L_g$ in the metric $g$, or the closely related Yamabe invariant. These functionals played a role in the solution of the Yamabe problem (see e.g. \cite{LP} and \cite{Sc}), and similar functionals for the Dirac operator have been studied in for instance \cite{Am1} and \cite{Am2}.

Note that in the paper \cite{PR} S. Paycha and S. Rosenberg have recently described a procedure for constructing conformal invariants via the canonical trace, and these invariants may also be studied by our methods presented in this paper.

Let us also mention that our main Theorem \ref{universal} has applications in the theory of conformally compact Einstein manifolds $(X^{n+1},g_+)$ with conformal infinities $M^n=\partial X$ (see e.g. \cite{FG1}). By \cite{GL} the $g_+$ are parametrized by the conformal infinities $(M^n,g)$, for metrics $g$ near the round metric $(S^n,g_0)$. It thus makes sense to study the variational problem for the renormalized volume, related to the gravitational action in AdS/CFT, which for $n$ odd is the constant term $V$ in the volume expansion
\[
\Vol_{g_+}(\{r>\varepsilon\})=c_0\varepsilon^{-n}+c_2\varepsilon^{-n+2}+\ldots+c_{n-1}\varepsilon^{-1}+V+o(1),\quad\textrm{as}\quad\varepsilon\to0,
\]
where $r$ denotes a boundary defining function for the conformal infinity $M^n$. Namely, for $n$ odd the renormalized volume $V$ is conformally invariant, and thus our rigidity result for the second variation holds, and given that the second variation is non-trivial, we obtain that the renormalized volume is locally extremalized at the round sphere. The first variation of $V$ for $n$ odd has been studied in \cite{Al} and \cite{An}, and for the analogous quantity in even dimension in \cite{CQY}.

It is also worth noticing, that even when a functional is not conformally invariant, but when we know explicitly how
it transforms under conformal changes of the metric (e.g. there is a Polyakov-type
formula), then we may in principle subtract this explicit term in order to obtain a conformally
invariant functional. This latter may then be analyzed by the results in this paper. Such methods might show to be interesting for instance in connection with AdS/CFT theory and renormalized volumes in the case of even-dimensional conformal infinities.

As an interesting aside, the paper \cite{SS} by P. Sarnak and A. Str\"o{}mbergsson, concerned with extremals of the zeta function of the Laplacian on functions on certain flat tori, contains ideas reminiscent of those in the present, in that invariance of a functional under a certain symmetry group $G$ plays a role in proving positive definiteness of the second variation of the functional. In their case the group $G$ was finite, while here we deal instead with the conformal group of the sphere, which is the semi-simple Lie group $\SO{n+1,1}$.

As a final remark we note that the representation theoretical part of our results may be
extended to other rank one semisimple Lie groups, which could give similar results in
parabolic geometries such as CR geometry (see e.g. \cite{FH} and \cite{JL}); here one would consider functionals on the
moduli space of CR structures. In this case the CR-manifolds for which one could a priori hope to study similar problems, would be odd-dimensional spheres, where the corresponding structure group would be $G=\SU{n,1}$. Some of the relevant representation theory of this group was adressed in the work by T. Branson, G. Olafson and B. \O{}rsted (e.g. \cite{BOO}, see also Baston-Eastwood \cite{BE}). Likewise it seems feasible to study conformal invariants in Lorentzian signature, or more generally pseudo-Riemannian spaces, with $S^p\times S^q$ as the main example, with conformal group $G=SO(p+1,q+1)$. In this connection it is worth pointing out the method of using analytic continuation in the signature of the metric in \cite{BO1} and \cite{BO3}, for extending the definition of conformal invariants to mixed signatures.

\bibliographystyle{amsalpha}

\end{document}